\documentclass{article}

\newtheorem{thm}{Theorem}
\newtheorem{cor}[thm]{Corollary}
\newtheorem{lem}[thm]{Lemma}

\begin{document}
\title{ Pairwise Relative Primality of Positive Integers}
\author{Jerry Hu}
\date{}
\maketitle

A classic result in number theory is that the probability that two positive integers are relatively prime is $6/\pi^2$.
More generally the probability that $k$ positive integers chosen arbitrarily and independently are relatively prime
is $ 1/\zeta(k)$, where $\zeta(k)$ is Riemann's zeta function. A short accessible proof of this result was given by J. E. Nymann \cite{nymann}. In a recent paper L. T\'{o}th \cite{toth} solved the problem of finding the probability that $k$ positive integers are pairwise relatively prime by the recursion method, he proved that, for $k \geq 2,$ the probability that $k$ positive integers are pairwise relatively prime is

$$ A_k = \prod_p \left(1-\frac{1}{p}\right)^{\hspace{-0.04in}k-1}\hspace{-0.05in}\left(1+\frac{k-1}{p}\right).$$

 Given a graph $G=(V,E)$ with $V=\{1,2,...,k\}$, the $k$ positive integers $a_1,a_2, ...,a_k$ are $G$-wise relatively prime if $(a_i, a_j)=1$ for $\{i,j\} \in  E$.  In this note we consider the problem of finding the probability $A_G$ that k positive integers are $G$-wise relatively prime.

For a $k$-tuple positive integers $u = (u_1,u_2,...,u_k)$,   let $ Q_G^{(u)}(n)$ denote the number of $k$-tuples of positive integers $a_1,a_2, ...,a_k $ with $1 \leq a_1,a_2,..., a_k \leq n$ such that $(a_i,u_i)=1$ for $i=1,...,k$ and they are $G$-wise relatively prime.

The next theorem gives an asymptotic formula for $Q_G^{(u)}(n)$ and the exact values of $A_G$. Before we state the theorem, let us introduce some notations. Given a graph $G=(V,E)$ with $V=\{1,2,...,k\}$, a subset $S \subset V$ is called independent if no two vertices of S are adjacent in $G$. We denote by $i_m(G)$ the number of independent sets of cardinality m in G, and for a subset $S \subset V$, we denote by $i_{m,S}(G)$ the number of independent sets of cardinality m in $G$ which contains at least one vertex in $S$. For a $k$-tuple positive integers $u = (u_1,u_2,...,u_k),$ and an integer $d$, the set of positive integers $i$ with $1\leq i \leq k$  such that $d$ divides $u_i$ is denoted by $S(u,d)$.

\begin{thm} For a graph $G=(V,E)$ with $V=\{1,2,...,k\}$, we have uniformly for $n, u_i \geq 1,  $
\begin{equation}  Q_G^{(u)}(n)=A_Gf_G(u)n^k + O(\theta(u)n^{k-1}\log^{k-1}n ),  \label{meqn}
\end{equation}
where
$$A_G=  \prod_p\left(\sum_{m=0}^{k}i_m(G)\left(1-\frac{1}{p}\right)^{\hspace{-0.04in}k-m}\hspace{-0.05in}\frac{1}{p^m}\right),$$
$$f_G(u)=\prod_{p|u_1u_2\cdot\cdot\cdot u_k}\left(1-\frac{\sum_{m=0}^{k}i_{m,S(u,p)}(G)(p-1)^{k-m}}{\sum_{m=0}^{k}i_m(G)(p-1)^{k-m}}\right),$$
and if $\theta(u_i)$ denotes the number of square free divisors of $u_i$, then $\theta (u)=\max\{ \theta(u_i), i=1,2,...,k \}.$
\end{thm}

\begin{cor}
The probability that $k$ positive integers $a_1, a_2, ..., a_k$ are $G$-wise relatively prime and $(a_i,u_i)=1$ for  $i=1,...,k$  is
$$ \lim_{n \rightarrow \infty} \frac{Q_G^{(u)}(n)}{n^k}=A_G f_G(u).$$
For $u_i=1$, the probability that $k$ positive integers are $G$-wise relatively prime is
$$  A_G=  \prod_p\left(\sum_{m=0}^{k}i_m(G)\left(1-\frac{1}{p}\right)^{\hspace{-0.04in}k-m}\hspace{-0.05in}\frac{1}{p^m}\right).$$
\end{cor}

In \cite{moree}, P. Moree proposed the problems of finding probabilities that k positive integers have exact (or at least) r relatively prime pairs. Using Theorem 1 and an Inclusion-Exclusion formula in combinatorics (see exercise 1 of chapter 2 in \cite{stanley}), we can give a solution to his problems.
\begin{cor}
The probability that $k$ positive integers have exactly $r$ relatively prime pairs  is
$$ A_{k,=r} = \sum_{i=r}^{k(k-1)/2} (-1)^{i-r} \left(\hspace{-0.05in}\begin{array}{c} i \\ r  \end{array}\hspace{-0.05in}\right)B_{k,i},$$
and the probability that $k$ positive integers have at least $r$ relatively prime pairs  is
$$ A_{k,\geq r} =  \sum_{i=r}^{k(k-1)/2} (-1)^{i-r} \left(\hspace{-0.05in}\begin{array}{c} i-1 \\ r-1  \end{array}\hspace{-0.05in}\right)B_{k,i},$$
where
$$ B_{k,i}= \sum_{|E| = i } A_G.$$
In particular, the probability that $k$ positive integers are pairwise not relatively prime is
$$ A_{k,=0}= \sum_{i=0}^{k(k-1)/2} (-1)^{i} B_{k,i}.$$
\end{cor}

In \cite{fernandez}, J. L. Fern\'{a}ndez and P. Fern\'{a}ndez proved that the number of relatively prime pairs is asymptotically normal as $k$ tends to $\infty$.

To prove Theorem 1 we need the following lemmas.

\begin{lem}
For $k,n \geq 1$, a graph $G=(V,E)$ with $V=\{1,2,...,k+1\}$, and $u=(u_1,u_2,...,u_{k+1})$ with $u_i \geq 1$,
$$ Q_G^{(u)}(n)  =  \sum_{\scriptsize \begin{array}{c} j=1\\(j,u_{k+1})=1 \end{array}}^nQ_{G-v}^{(j*u)}(n),$$
where  $G-v$ is the graph obtained from G by deleting the vertex$\mbox{ v=k+1}$ together with all the edges incident to $v$, and if $(j*u)_i$ denotes the $i$th component of $j*u$, then
$$ (j*u)_i= \left\{ \begin{array}{ll} ju_i & \mbox{{\rm  if} i {\rm is adjacent to }} v \mbox{ {\rm in }} G, \\ u_i & {\rm otherwise},\end{array} \right. $$
for $i=1,2,...,k.$
\end{lem}
{\em Proof.} The $k+1$ positive integers $a_1,a_2,...,a_{k+1}$ are $G$-wise relatively prime and $ (a_i,u_i)=1$ for  $i=1,2,...,k+1$ if and only if the first $k$ positive integers $a_1,a_2,...,a_k$ are $(G-v)$-wise relatively prime and $(a_i,u_i)=1$ for $i=1,2,...,k,$ and $(a_i,a_{k+1})=1$ when the vertex $i$ is adjacent to the vertex $ v=k+1,$ and $(a_{k+1}, u_{k+1})=1$. We have
$$ Q_G^{(u)}(n)= \sum_{\scriptsize \begin{array}{c} a_{k+1}=1\\(a_{k+1},u_{u+1})=1 \end{array}}^nQ_{G-v}^{(a_{k+1}*u)}(n)=\sum_{\scriptsize \begin{array}{c} j=1\\(j,u_{k+1})=1 \end{array}}^nQ_{G-v}^{(j*u)}(n).$$

\begin{lem}
For $k, u_i \geq 1$, a graph $G=(V,E)$ with $V=\{1,2,...,k\}$, and $S$ a subset of vertices in $V$,
$$ \frac{f_G(j*u)}{f_G(u)}= \sum_{d|j}\frac{\mu(d)\alpha_{G,S}(u,d)}{\alpha_G(u,d)},$$
if $d$ is square free, then
$$ \frac{\alpha_{G,S}(u,d)}{\alpha_G(u,d)} \leq \frac{k^{\omega(d)}}{d},$$
where
$$ \alpha_G(u,d) = \prod_{p\,| d}\left(\sum_{m=0}^{k}i_m(G-S(u,p))(p-1)^{k-m}\right),$$
$$\alpha_{G,S}(u,d)=\prod_{p\,|d}\left(\sum_{m=0}^{k}i_{m,S}(G-S(u,p))(p-1)^{k-m}\right),$$
and if $(j*u)_i$ denotes the $i$th component of $j*u$, then
$$ (j*u)_i= \left\{ \begin{array}{ll} ju_i & \mbox{{\rm  if} i {\rm is in }} S, \\ u_i & {\rm otherwise},\end{array} \right. $$
for $i=1,2,...,k,$  and $\omega(d)$ denote the number of distinct prime factors of $d.$
\end{lem}
{\em Proof.} It suffices to  verify the equality for $j=p^a$ a prime power:

\begin{eqnarray*}
\lefteqn{\sum_{d | p^a}\frac{\mu(d)\alpha_{G,S}(u,d)}{\alpha_G(u,d)}}\\&  = & 1-\frac{\sum_{m=0}^{k}i_{m,S}(G\hspace{-0.03in}-\hspace{-0.03in}S(u,p))(p-1)^{k-m}}{\sum_{m=0}^{k}i_m(G\hspace{-0.03in}-\hspace{-0.03in}S(u,p))(p-1)^{k-m}}
 \\
& =&  \frac{\sum_{m=0}^{k}i_m(G\hspace{-0.03in}-\hspace{-0.03in}S(u,p))(p-1)^{k-m}\hspace{-0.02in}-\hspace{-0.03in}\sum_{m=0}^{k}i_{m,S}(G\hspace{-0.03in}-\hspace{-0.03in}S(u,p))(p-1)^{k-m}}{\sum_{m=0}^{k}i_m(G-S(u,p))(p-1)^{k-m}}\\
&=&\frac{\sum_{m=0}^{k}i_m(G)(p-1)^{k-m}-\sum_{m=0}^{k}i_{m,S\cup S(u,p)}(G)(p-1)^{k-m}}{\sum_{m=0}^{k}i_m(G)(p-1)^{k-m}-\sum_{m=0}^{k}i_{m,S(u,p)}(G)(p-1)^{k-m}}\\
&=& \frac{f_{G}(p^a*u)}{f_G(u)}.
\end{eqnarray*}

Now we prove the inequality. Notice that
$$i_0(G-S(u,p))=1, \hspace{0.3in} i_{0,S}(G-S(u,p))=0, \hspace{0.3in} i_{m,S}(G-S(u,p)) \leq i_m(G-S(u,p)).$$
Then when $d$ is square free,
\begin{eqnarray*}
\frac{\alpha_{G,S}(u,d)}{\alpha_G(u,d)} & = & \prod_{p\,| d}\frac{\sum_{m=0}^{k}i_{m,S}(G-S(u,p))(p-1)^{k-m}}{\sum_{m=0}^{k}i_m(G-S(u,p))(p-1)^{k-m}}\\
& \leq &  \prod_{p\,| d}\frac{\sum_{m=1}^{k}i_m(G-S(u,p))(p-1)^{k-m}}{(p-1)^k + \sum_{m=1}^{k}i_m(G-S(u,p))(p-1)^{k-m}}\\
& \leq &  \prod_{p\,| d}\frac{\sum_{m=1}^{k}\left(\hspace{-0.05in}\begin{array}{c} k \\ m  \end{array}\hspace{-0.05in}\right)(p-1)^{k-m}}{(p-1)^k+\sum_{m=1}^{k}\left(\hspace{-0.05in}\begin{array}{c} k \\ m  \end{array}\hspace{-0.05in}\right)(p-1)^{k-m}}\\
& =& \prod_{p\,| d}\frac{(p^k-(p-1)^k)}{p^k}\\
&\leq & \prod_{p\,| d}\frac{kp^{k-1}}{p^k}\\
&=& \prod_{p\,| d}\frac{k}{p}= \frac{k^{\omega(d)}}{d}\\
\end{eqnarray*}

For the proof of the theorem, we proceed by induction on $k$. For $k=1, $ we have by the Inclusion-Exclusion Principle
\begin{eqnarray*}
Q_{\{1\}}^{(u_1)}(n) & = & \sum_{\scriptsize \begin{array}{c} j=1 \\(j,u_1)=1 \end{array}}^{n} 1=\sum_{d|u_1}\mu(d)\lfloor\frac{n}{d}\rfloor=\sum_{d|u_1}\mu(d)(\frac{n}{d}+O(1))\\
 & = & n\sum_{d|u_1}\frac{\mu(d)}{d}+O(\sum_{d|v}\mu^2(d)).
\end{eqnarray*}
Hence,
\begin{equation} Q_{\{1\}}^{(u_1)}(n)  =  \sum_{\scriptsize \begin{array}{c} j=1 \\(j,u_1)=1 \end{array}}^{n}1=n\frac{\phi(u_1)}{u_1}+O(\theta(u_1)) \label{eueqn} \end{equation}
and (\ref{meqn}) is true for $k=1$ with $A_{\{1\}}=1, f_{\{1\}}(u_1)=\frac{\phi(u_1)}{u_1}, \phi$ denoting the Euler function.

Suppose that (\ref{meqn}) is valid for $k$, we prove it for $k+1$. Let $u=(u_1,u_2,...,u_{k+1})$ and $u'=(u_1,u_2,...,u_k)$, from Lemma $4$ we have

\begin{eqnarray}
Q_G^{(u)}(n) & = & \sum_{\scriptsize \begin{array}{c} j=1\\(j,u_{k+1})=1 \end{array}}^nQ_{G-v}^{(j*u)}(n)\nonumber \\
& = & \sum_{\scriptsize \begin{array}{c} j=1\\(j,u_{k+1})=1 \end{array}}^n A_{G-v} f_{G-v}(j*u)n^k+O(\theta(j*u)n^{k-1}\log^{k-1}\hspace{-0.05in}n) \nonumber  \\
 & = & A_{G-v} f_{G-v}(u')n^k \sum_{\scriptsize \begin{array}{c} j=1\\(j,u_{k+1})=1 \end{array}}^n \frac{f_{G-v}(j*u')}{f_{G-v}(u')}  \nonumber \\
& & \hspace{1.1in} \mbox{} + O(\theta(u')n^{k-1}\log^{k-1}\hspace{-0.05in}n\sum_{j=1}^n\theta(j)).\label{eqnb}
\end{eqnarray}
Here $\sum_{j=1}^n\theta(j) \leq \sum_{j=1}^n\tau_2(j)=O(n\log n),$ where $\tau_2 = \tau$ is the divisor \mbox{function}.

Furthermore, in Lemma $5$ choosing the subset $S$ to be the open neighbourhood $N(v)$ of $v$, which is the set of vertices adjacent to $v$,  we have
\begin{eqnarray*}
\lefteqn{ \sum_{\scriptsize \begin{array}{c} j=1\\(j,u_{k+1})=1 \end{array}}^n  \frac{f_{G-v}(j*u')}{f_{G-v}(u')} }  \\
      & = & \sum_{\scriptsize \begin{array}{c} de=j \leq n\\(j,u_{k+1})=1 \end{array}}\frac{\mu(d)\alpha_{G-v,N(v)}(u',d)}{\alpha_{G-v}(u',d)} \\
  &  = & \sum_{\scriptsize \begin{array}{c} d\leq n\\ (d,u_{k+1})=1 \end{array}}
 \frac{\mu(d)\alpha_{G-v,N(v)}(u',d)}{\alpha_{G-v}(u',d)}
\sum_{\scriptsize \begin{array}{c} e\leq\frac{n}{d}\\(e, u_{k+1})=1\end{array}}\hspace{-0,1in}1
\end{eqnarray*}

Using  (\ref{eueqn}), we have
\begin{eqnarray}
\lefteqn{ \sum_{\scriptsize \begin{array}{c} j=1\\(j,u_{k+1})=1 \end{array}}^n \hspace{-0.1in} \frac{f_{G-v}(j*u')}{f_{G-v}(u')}}  \nonumber \\
  & =&\hspace{-0.2in} \sum_{\scriptsize \begin{array}{c} d\leq n\\ (d,u_{k+1})=1 \end{array}}
  \hspace{-0.2in} \frac{\mu(d)\alpha_{G-v,N(v)}(u,d)}{\alpha_{G-v}(u',d)}
 \left(\frac{\phi(u_{k+1})}{u_{k+1}}\frac{n}{d}+O(\theta(u_{k+1}))\right) \nonumber \\
 & =& \frac{\phi(u_{k+1})}{u_{k+1}} n \hspace{-0.2in}
 \sum_{\scriptsize \begin{array}{c} d\leq n\\ (d,u_{k+1})=1 \end{array}} \hspace{-0.2in}
 \frac{\mu(d)\alpha_{G-v,N(v)}(u',d)}{d \alpha_{G-v}(u',d)}
 +   O\left(\theta(u_{k+1})\sum_{d\leq n}\frac{k^{\omega(d)}}{d}\right), \label{eqna}
\end{eqnarray}
by Lemma $5.$

Hence, the main term  of (\ref{eqna}) is
\begin{eqnarray*}
\lefteqn{\frac{\phi(u_{k+1})}{u_{k+1}}n \sum_{\scriptsize \begin{array}{c} d=1\\ (d,u_{k+1})=1 \end{array}}^{\infty}\frac{\mu(d)\alpha_{G-v,N(v)}(u',d)}{d \alpha_{G-v}(u',d)} } \\
 &=&\frac{\phi(u_{k+1})}{u_{k+1}}n\prod_{p\hspace{-0.01in} \not \hspace{0.025in}|  \hspace{0.01in}u_1u_2\cdot\cdot\cdot u_ku_{k+1}}
 \left(1-\frac{\sum_{m=0}^{k}i_{m,N(v)}(G-v)(p-1)^{k-m}}{p\,(\sum_{m=0}^{k}i_m(G-v)(p-1)^{k-m})}\right)\\
 & & \prod_{\scriptsize \begin{array}{c} p |u_1u_2\cdot\cdot\cdot u_k\\p \! \! \not \hspace{-0.008in} |  \, u_{k+1} \end{array}}  \left(1-\frac{\sum_{m=0}^{k}i_{m,N(v)}(G-v-S(u,p))(p-1)^{k-m}}{p\,(\sum_{m=0}^{k}i_m(G-v-S(u,p))(p-1)^{k-m})}\right)  \\
 &=& n\prod_{p}\left(1-\frac{\sum_{m=0}^{k}i_{m,N(v)}(G-v)(p-1)^{k-m}}{p\,(\sum_{m=0}^{k}i_m(G-v)(p-1)^{k-m})}\right)   \\
   & &\prod_{p \, | u_{k+1}}\left(1-\frac{1}{p}\right) \prod_{p \, |  u_1u_2\cdot\cdot\cdot u_{k+1}}
 \left(1-\frac{\sum_{m=0}^{k}i_{m,N(v)}(G-v)(p-1)^{k-m}}{p\,(\sum_{m=0}^{k}i_m(G-v)(p-1)^{k-m})}\right)^{-1} \\
 & & \prod_{\scriptsize \begin{array}{c} p |u_1u_2\cdot\cdot\cdot u_k\\p \! \! \not \hspace{-0.01in} |  \, u_{k+1} \end{array}}  \left(1-\frac{\sum_{m=0}^{k}i_{m,N(v)}(G-v-S(u,p))(p-1)^{k-m}}{p\,(\sum_{m=0}^{k}i_m(G-v-S(u,p))(p-1)^{k-m})}\right),
     \\
  \end{eqnarray*}

 and its O-terms are
 \begin{eqnarray*}
 O(n\sum_{d > n}\frac{k^{\omega(d)}}{d^2} ) &=&
      O(n\sum_{d > n}\frac{\tau_{k}(d)}{d^2}) \\
     & = & O(\log^{k-1}\hspace{-0.05in}n)
 \end{eqnarray*}
by Lemma 3(b) in \cite{toth}, which gives an asymptotic estimate of the sum $$ \sum_{n > x}\frac{\tau_k(n)}{n^2}=O(\frac{\log^{k-1}x}{x})$$
  and

\begin{eqnarray*}
 O(\theta(u_{k+1})\sum_{d \leq n}\frac{k^{\omega(d)}}{d})
& = & O(\theta(u_{k+1})\sum_{d\leq n}\frac{\tau_{k}(d)}{d})\\
&=& O(\theta(u_{k+1})\log^{k}\hspace{-0.02in}n)
\end{eqnarray*}
from Lemma 3(a) in \cite{toth}, which gives an asymptotic estimate of the sum $$ \sum_{n\leq x}\frac{\tau_k(n)}{n}=O(\log^{k}\hspace{-0.02in}x). $$

Substituting into (\ref{eqnb}), we get
\begin{eqnarray*}
Q_{G}^{(u)}(n) & = & A_{G-v}\prod_{p}\left(1-\frac{\sum_{m=0}^{k}i_{m,N(v)}(G-v)(p-1)^{k-m}}{p\,(\sum_{m=0}^{k}i_m(G-v)(p-1)^{k-m})}\right)   \\
 & & f_{G-v}(u)\prod_{p \, | u_{k+1}}\left(1-\frac{1}{p}\right) \prod_{p \, |  u_1u_2\cdot\cdot\cdot u_{k+1}}
 \left(1-\frac{\sum_{m=0}^{k}i_{m,N(v)}(G-v)(p-1)^{k-m}}{p\,(\sum_{m=0}^{k}i_m(G-v)(p-1)^{k-m})}\right)^{\hspace{-0.05in}-1} \\
 & & \prod_{\scriptsize \begin{array}{c} p |u_1u_2\cdot\cdot\cdot u_k\\p \! \! \not \hspace{-0.01in}|  \, u_{k+1} \end{array}}  \left(1-\frac{\sum_{m=0}^{k}i_{m,N(v)}(G-v-S(u,p))(p-1)^{k-m}}{p\,(\sum_{m=0}^{k}i_m(G-v-S(u,p))(p-1)^{k-m})}\right) n^{k+1},
     \\
  & & \mbox{}+O(n^k\log^{k-1}\hspace{-0.05in}n)+O(\theta(u_{k+1})n^k\log^{k}\hspace{-0.05in}n)+O(\theta(u')n^k\log^{k}\hspace{-0.05in}n)\\
 &= & A_{G}f_{G}(u)n^{k+1}+O(\theta(u)n^k\log^{k}\hspace{-0.05in}n)
\end{eqnarray*}
by a simple computation, which shows that the formula is true for $k+1$ and we complete the proof.

T. Freiberg computed the probability that three positive integers are pairwise not relatively prime, see \cite{moree}. As an application of our method, we compute the probability that four positive integers are pairwise not relatively prime.

Let $A_{3,i}$ denote the probability that three positive integers have $i$ relatively prime pairs, for $i=1,2,3.$
By Theorem $1$ and Corollary $3$, we have
$$A_{3,0}=1, \hspace{0.1in}A_{3,1}=\prod_p \left(1-\frac{1}{p^2}\right),\hspace{0.1in} A_{3,2}=\prod_p \left(1-\frac{1}{p}\right)\left(1+\left(1-\frac{1}{p}\right)\frac{1}{p}\right),$$
$$ \hspace{-2.6in}A_{3,3} = \prod_p  \left(1-\frac{1}{p}\right)^2\left(1-\frac{2}{p}\right), $$
and
$$ A_{3, =0}= 1- 3A_{3,1}+3A_{3,2}-A_{3,3}.$$
This recovers T. Freiberg's result.

When $k=4$, the number of graphs with given number of edges and the number of independent sets are summarized in the following table.

$$\begin{tabular}{|l|l|l|} \hline
number of edges & number of graphs & number of independent sets \\ \hline
0 & 1 & $i_2(G)=6, i_3(G)=4, i_4(G)=1$ \\ \hline
1&  6 & $i_2(G)=6, i_3(G)=4, i_4(G)=1$ \\ \hline
2 & type I: 12 & type I: $ i_2(G)=4, i_3(G)=2, i_4(G)=0$ \\ \cline{2-3}
  &  type II: 3  & type II: $i_2(G)=4, i_3(G)=0, i_4(G)=0$ \\ \hline
3 & type I: 4 & type I: $ i_2(G)=3, i_3(G)=1, i_4(G)=0$ \\ \cline{2-3}
  &  type II: 16  & type II: $i_2(G)=3, i_3(G)=0, i_4(G)=0$ \\ \hline
4 & 15 & $ i_2(G)=2, i_3(G)=0, i_4(G)=0$ \\ \hline
5 & 6 & $i_2(G)=1, i_3(G)=0, i_4(G)=0$ \\ \hline
6 & 1 & $i_2(G)=0, i_3(G)=0, i_4(G)=0$ \\ \hline
\end{tabular} $$

Given graph $G$ with four vertices and $i$ edges, let $A_{4,i}$ denote the probability that integers are $G$-wise relatively prime for $i=1,2,...,6,$ if there are more than one type of graphs with fixed number of edges $i$, for the graphs of type $j$, the probability is denoted by $A_{4,i,j}$. Again by Theorem $1$ and Corollary $3$, we have
$$A_{4,0}=1, A_{4,1}=\prod_p \left(1-\frac{1}{p^2}\right), \hspace{0.08in}A_{4,2,1}=\prod_p \left(1-\frac{1}{p}\right)\left(1+\left(1-\frac{1}{p}\right)\frac{1}{p}\right),$$ 
$$A_{4,2,2}=\prod_p \left(1-\frac{1}{p^2}\right)^2, \hspace{0.1in}
 A_{4,3,1}=\prod_p \left(1-\frac{1}{p}\right)\left(1+\left(1-\frac{1}{p}\right)^2 \frac{1}{p}\right), $$
 $$ A_{4,3,2}=\prod_p  \left(1-\frac{1}{p}\right)^2\left(1-\frac{2}{p}\right),\hspace{0.08in} A_{4,4}=\prod_p  \left(1-\frac{1}{p}\right)^2\left(1+\frac{2}{p}-\frac{1}{p^2}\right),$$
$$A_{4,5}=\prod_p  \left(1-\frac{1}{p}\right)^2\left(1+ 2\left(1-\frac{1}{p}\right)\frac{1}{p}\right), \hspace{0.08in}A_{4,6}=\prod_p  \left(1-\frac{1}{p}\right)^3\left(1+\frac{3}{p}\right), $$ and
 $$A_{4,=0}=1-6A_{4,1}+12A_{4,2,1}+3A_{4,2,2}-4A_{4,3,1}-16A_{4,3,2}+15A_{4,4}-6A_{4,5}+A_{4,6}.$$

{\footnotesize Department of Mathematics and Computer Science,
 School of Arts \& Sciences, University of Houston-Victoria,  14000 University Blvd, Sugar Land, TX 77479, USA}

\end{document}